# Extended Lorenz majorization and frequencies of distances in an undirected network

Leo Egghe

Hasselt University, Belgium

E-mail: leo.egghe@uhasselt.be

**Abstract**

Purpose – To contribute to the study of networks and graphs.

Design/methodology/approach – We apply standard mathematical thinking.

Findings – We show that the distance distribution in an undirected network Lorenz majorizes the one of a chain. As a consequence, the average and median distances in any such network are smaller than or equal to those of a chain.

Research limitations – We restricted our investigations to undirected, unweighted networks.

Practical implications – We are convinced that these results are useful in the study of small worlds and the so-called six degrees of separation property.

Originality/value – To the best of our knowledge our research contains new network results, especially those related to frequencies of distances.

## Introduction

Let G=(*V,E*) be an undirected network, where $V = (v_k)_{k=1,\ldots,N}$ denotes the set of nodes or vertices and *E* denotes the set of links or edges. As collaboration (of scientists, universities, countries), bibliographic coupling (of articles, books), and co-citation (of articles, books) are all examples of undirected networks, it goes without saying that the study of these networks is of great importance for bibliometrics (Rousseau et al., 2018).

We assume that #V = N > 1. A chain in a network is a sequence of different nodes one by one connected by edges. The distance d between two nodes is equal to the number of links situated on a shortest chain (often called a shortest path) between these two nodes. Consequently, the distance between two nodes connected by an edge is equal to one. Each network studied in this article is assumed to be connected, i.e. there is a chain between any two nodes. Hence, for each node, there exists another node at a distance one. The total number of distances between any pair of nodes in this network is equal to N(N-1)/2, where, for $v_1, v_2 \in V$, $d(v_1,v_2)=d(v_2,v_1)$ is considered only once.

Notation. We denote by $\alpha_j$, j= 1,…, N-1, the number of times distance j occurs in network G. The array $A = (\alpha_1, \alpha_2, \ldots, \alpha_{N-1})$ is called the $\alpha -$ array of the network G.

Some immediate properties

1) $\sum_{j=1}^{N-1} \alpha_j = \frac{N(N-1)}{2}$  (1)

2) $\alpha_1 = \#E \geq N - 1$  (2)

3) $\alpha_{N-1} = 0 \text{ or } 1$  (3)

4) $\alpha_j = 0 \Rightarrow \forall k \geq j: \alpha_k = 0$  (4)

We give a short proof of

5) $\alpha_2 \leq \frac{(N-1)(N-2)}{2}$ (5)

Indeed, if $\alpha_2 > \frac{(N-1)(N-2)}{2}$ then $\frac{N(N-1)}{2} = \sum_{j=1}^{N-1} \alpha_j \geq \alpha_1 + \alpha_2 > (N-1) + \frac{(N-1)(N-2)}{2} = \frac{N(N-1)}{2}$, which is a contradiction.

6) In a chain of length j there are 2 chains of length j-1, 3 chains of length j-2, and k chains of length j-k+1 (0 < k < j).

7) The distance frequency array of a complete N-node network, K_N, is $\left(\frac{N(N-1)}{2}, 0, \ldots, 0\right)$.

We recall the definition of the majorization order (Hardy et al., 1934). Let $X = (x_j)$ and $Y = (y_j)$, j=1,…, N-1 be two (N-1)-sequences of non-negative numbers, ordered decreasingly then X majorizes Y, denoted as $X \succcurlyeq Y$, if

$$\forall i = 1, \ldots, N-2: \sum_{j=1}^{i} x_j \geq \sum_{j=1}^{i} y_j \qquad (6)$$

and $\qquad \sum_{j=1}^{N-1} x_j = \sum_{j=1}^{N-1} y_j \qquad (7)$

We recall that if $X \succcurlyeq Y$ then the (standard) Lorenz curve of X (Lorenz, 1905) is situated above the Lorenz curve of Y. We now extend the notions of majorization and Lorenz curve by removing the requirement to be arranged in decreasing order.

Definition. Extended majorization order

Let $X = (x_j)$ and $Y = (y_j)$, j=1,…, N-1 be two (N-1)-sequences of non-negative numbers, then X majorizes Y (in the extended sense), denoted as $X \succcurlyeq Y$ (we keep the same notation), if

$$\forall\, i = 1, \ldots, N-2: \sum_{j=1}^{i} x_j \geq \sum_{j=1}^{i} y_j \qquad (8)$$

and
$$\sum_{j=1}^{N-1} x_j = \sum_{j=1}^{N-1} y_j \qquad (9)$$

Definition. Extended Lorenz curve

Let $X = (x_j)$ be an (N-1)-sequence of non-negative numbers and let $s_i = \sum_{j=1}^{i} x_j$ be the j$^{th}$ partial sum. Hence $s_{N-1}$ = TOT denotes the total sum of the numbers in X; $s_0$ is set equal to 0. Now plot the points $\left(\frac{i}{N-1}, \frac{s_i}{TOT}\right)_{i=0,\ldots,N-1}$ and connect them by line segments to obtain a curve joining the origin (0,0) with the point (1,1). We refer to this curve as the extended Lorenz curve. Contrary to the standard Lorenz curve this curve is not necessarily concave (but of course still increasing). An example is shown in Fig. 3. If X is increasing then the extended Lorenz curve coincides with the classical Lorenz curve. If $X \succcurlyeq Y$ then the extended Lorenz curve of X is situated above the extended Lorenz curve of Y.

**The main result**

Given the number of nodes, N, we next show a majorization result between the frequency sequence of a complete network $K_N$, that of a general network G=(V, E), denoted as *A*, and the frequency sequence *C* of a chain.

Theorem

Given a network G=(V,E) with N nodes, then,

$$\left(\frac{N(N-1)}{2}, 0, \ldots, 0\right) \succcurlyeq A = (\alpha_1, \alpha_2, \ldots, \alpha_{N-1}) \succcurlyeq C = (N-1, N-2, \ldots, 1).$$

By (8) and (9) the second inequality means that

$$\forall\, i = 1, \ldots, N-2:\ \sum_{j=1}^{i} \alpha_j \geq \sum_{j=1}^{i}(N-j) \tag{10}$$

and

$$\sum_{j=1}^{N-1} \alpha_j = \sum_{j=1}^{N-1}(N-j) \tag{11}$$

We moreover prove that

$$\alpha_i \leq \sum_{j=i}^{N-1} \alpha_j \leq \frac{(N-i+1)(N-i)}{2} =: \beta_i \tag{12}$$

Proof. By (1), we already know that

$$\frac{N(N-1)}{2} = \sum_{j=1}^{N-1} \alpha_j = \sum_{j=1}^{N-1}(N-j)$$

which proves (11). Assume now that $\alpha_i = 0$, for i > 1, then we already know that $\forall\, k \geq i: \alpha_k = 0$ and thus

$$\sum_{j=1}^{i} \alpha_j = \sum_{j=1}^{N-1} \alpha_j = \frac{N(N-1)}{2} \geq \sum_{j=1}^{i}(N-j)$$

Assume now that for some i > 1, $\alpha_i \neq 0$, then we will prove that also in this case $\sum_{j=1}^{i} \alpha_j \geq \sum_{j=1}^{i}(N-j)$. This will be done in several steps. First, we show that $\sum_{j=1}^{i} \alpha_j \geq \frac{i(i+1)}{2}$. Indeed: there exists in V at least one chain of length i, connecting nodes to which we refer as $u_1, u_2, \ldots, u_{i+1}$. Then, by property 6, we know that $\forall\, j = 1, \ldots, i: \alpha_j \geq i - j + 1$, and hence $\sum_{j=1}^{i} \alpha_j \geq \left(\sum_{j=1}^{i}(i - j + 1)\right) = \frac{i(i+1)}{2}$. In the next step, we show that this inequality can be refined to $\sum_{j=1}^{i} \alpha_j \geq \left(\frac{i(i+1)}{2}\right) + i$. Indeed, as the network under study is connected, there exists a node in the network, denoted as $u_{i+2}$, connected to at least one node of the chain $u_1, u_2, \ldots, u_{i+1}$. This point $u_{i+2}$ has a distance d, 0 < d ≤ i to at least i points in the chain. Now, adding the i distances involving the point $u_{i+2}$ we obtain $\sum_{j=1}^{i} \alpha_j \geq$ $\left(\frac{i(i+1)}{2}\right) + i$. We further note that $\forall j \leq i$, each point in the set S = { $u_1$, $u_2$, ..., $u_{i+2}$} has at least j points in the set S at a distance $0 < d \leq j$.

Now we continue in this way. Assuming that we have a set T of i+n connected nodes {$u_1$, ..., $u_{i+1}$, $u_{i+2}$, ..., $u_{i+n+1}$} from which we already derived that $\sum_{j=1}^{i} \alpha_j \geq \left(\frac{i(i+1)}{2}\right) + (ni)$ and for which we know that $\forall j \leq i$ each point in the set T, has at least j points in the set T at a distance $0 < d \leq j$. We again apply connectedness to get a new node $u_{i+n+2}$ at a distance d, $0 < d \leq i$ to all points in T, leading to $\sum_{j=1}^{i} \alpha_j \geq \left(\frac{i(i+1)}{2}\right) + (n+1)i$. Again we observe that $\forall j \leq i$, each point in T* = {$u_1$, ..., $u_{i+1}$, $u_{i+2}$, ..., $u_{i+n+2}$ } has a distance d, $0 < d \leq j$, with at least j points in the set T*. This procedure ends with n = N-i-2 for which $\sum_{j=1}^{i} \alpha_j \geq \left(\frac{i(i+1)}{2}\right) + (N-i-1)i = Ni - \sum_{j=1}^{i} j = \sum_{j=1}^{i}(N-j)$ which proves the inequality in the case $\alpha_i \neq 0$, and hence (10).

Now we prove (12). Using (10) and (11) we have:

$$\forall i = 1, \ldots, N-1: \alpha_i \leq \sum_{j=i}^{N-1} \alpha_j = \sum_{j=1}^{N-1} \alpha_j - \sum_{j=1}^{i-1} \alpha_j \leq \frac{N(N-1)}{2} - \sum_{j=1}^{i-1}(N-j) = \beta_i \qquad (12)$$

where we still have to prove the final inequality in (12). For this, we first observe that:

$$\beta_i - \beta_{i+1} = \frac{(N-i+1)(N-i)}{2} - \frac{(N-i)(N-i-1)}{2} = (N-i)$$

Now,

$$\frac{N(N-1)}{2} - \sum_{j=1}^{i-1}(N-j) = \frac{N(N-1)}{2} - \sum_{j=1}^{i-1}(\beta_j - \beta_{j+1})$$

$$= \frac{N(N-1)}{2} - (\beta_1 - \beta_i) = \beta_i \quad \square$$

Remarks and consequences

1) It follows immediately from the previous theorem that for a given number N and $\alpha$ −array $A = (\alpha_1, \alpha_2, \ldots, \alpha_{N-1})$ the median of A is smaller than or equal to the median of the chain of length N-1 (N nodes).

2) It is always possible to find a network with N (N > 3) nodes such that $\alpha_i < \beta_i$, and this for each number i=1,…, N-1. Indeed, consider for N > 3, a network for which $\forall i, i = 1, \ldots, N-3, \alpha_i \neq 0$; $\alpha_{N-2} = 1$ and $\alpha_{N-1} = 0$. In this case $\alpha_{N-1} = 0 < \beta_{N-1} = 1$; $\alpha_{N-2} = 1 < \beta_{N-2} = 3$ and, using (12), $\forall i, i = 1, \ldots, N-3, \alpha_i < \sum_{j=i}^{N-2} \alpha_j = \sum_{j=i}^{N-1} \alpha_j \leq \beta_i$. Such a network may look like shown in Fig. 1.

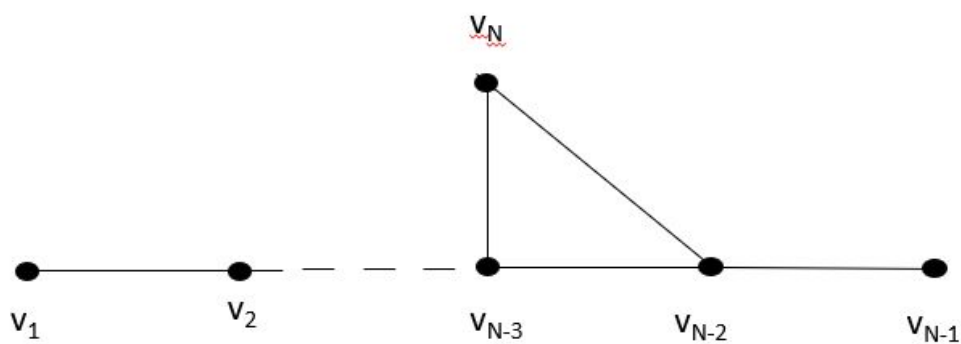

Fig. 1. An example of a network for which $\forall i = 1, \ldots, N-1: \alpha_i < \beta_i$.

3) The inequality $\alpha_i \leq \beta_i$ cannot be made more precise as for a chain of length i, $\alpha_i = \beta_i = 1$, for each i=1,…, N-1.

4) If N > 2, then it is impossible that $\forall\, i = 1, \ldots, N-1, \alpha_i = \beta_i$. Indeed: $\beta_1 = \frac{N(N-1)}{2}$ and if $\alpha_1 = \frac{N(N-1)}{2}$ then automatically $\alpha_i = 0, i = 2, \ldots, N-1$, while this is not the case for $\beta_i$.

5) If A is an array of length N-1, consisting of non-negative natural numbers such that

$$\left(\frac{N(N-1)}{2}, 0, \ldots, 0\right) \succcurlyeq A = (\alpha_1, \alpha_2, \ldots, \alpha_{N-1}) \succcurlyeq C = (N-1, N-2, \ldots, 1)$$

then the components of A do not have to be frequencies of distances in a network. Indeed, let N = 4 and let A = (4,1,1), then $(4,1,1) \succcurlyeq C = (3,2,1)$. Yet, there does not exist a network with (4,1,1) as distance frequencies: the third component is equal to one indicating that the network must be a chain but for a chain with 4 nodes, $\alpha_1 = 3$ and not 4.

Even if the last component of A is zero a counterexample is possible. Indeed, with N=5, we have $(4,3,3,0) \succcurlyeq C = (4,3,2,1)$. Such a network must have at least one chain of length three (connecting four nodes). The fifth node must be connected to the second or the third node in the chain. Hence A must necessarily be (4,4,2,0) and cannot be (4,3,3,0). These examples lead to the open question of finding the conditions under which such an array A is the frequency array of the distances in a (connected) network.

6) From the above and the main theorem we see that max {Md; Md is the median distance in an N-node network} is strictly smaller than max { $\bar{d}$: $\bar{d}$ is the average distance of an N-node network}. Although Md $\leq \bar{d}$ is not always true: a star with a center and N-1 rays (N>4) is an example (Md = 2 and $\bar{d} = \frac{2(N-1)}{N} < 2$ ), we have that if the $\alpha-$ sequence of a network is decreasing then clearly Md $\leq \bar{d}$. The reverse of this result does not hold in general. This is illustrated by G₀ (N=7) in Fig. 2 below. Its $\alpha$–sequence is not decreasing, namely (6,7,6,2,0,0) but yet Md = 2 < $\bar{d}$ = 46/21 ≈ 2.19.

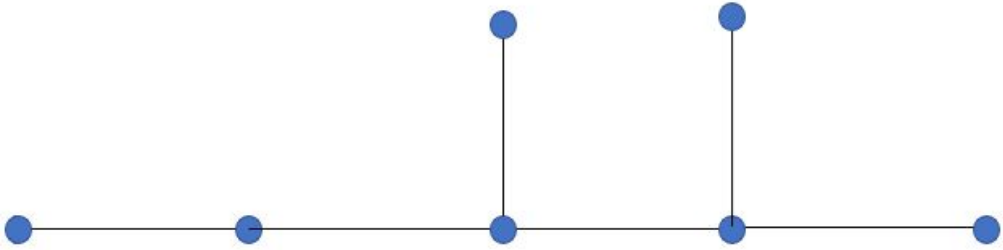

Fig. 2 An example of a network, G₀, with seven nodes

**A result about the average distance in a network**

If N is fixed and the array $A = (\alpha_1, \alpha_2, \ldots, \alpha_{N-1})$ denotes the frequencies of the distances in a network, then $\frac{2}{N(N-1)}\sum_{i=1}^{N-1} i\alpha_i$ denotes the average distance between nodes in this network, say $\bar{d}$.

Theorem. If $A^{(1)} = \left(\alpha_1^{(1)}, \alpha_2^{(1)}, \ldots, \alpha_{N-1}^{(1)}\right) \succcurlyeq A^{(2)} = \left(\alpha_1^{(2)}, \alpha_2^{(2)}, \ldots, \alpha_{N-1}^{(2)}\right)$ then $\overline{d_1} \leq \overline{d_2}$.

Proof. As $A^{(1)} = \left(\alpha_1^{(1)}, \alpha_2^{(1)}, \ldots, \alpha_{N-1}^{(1)}\right) \succcurlyeq A^{(2)} = \left(\alpha_1^{(2)}, \alpha_2^{(2)}, \ldots, \alpha_{N-1}^{(2)}\right)$, we know that $\forall\, i = 1, \ldots, N-2: \sum_{j=1}^{i} \alpha_j^{(1)} \geq \sum_{j=1}^{i} \alpha_j^{(2)}$ and $\sum_{j=1}^{N-1} \alpha_j^{(1)} = \sum_{j=1}^{N-1} \alpha_j^{(2)}$. Consequently: $\forall\, i = 2, \ldots, N-2: \sum_{j=i}^{N-1} \alpha_j^{(1)} \leq \sum_{j=i}^{N-1} \alpha_j^{(2)}$.

Now, $\overline{d_1} = \frac{2}{N(N-1)} \sum_{j=1}^{N-1} j\alpha_j^{(1)}$

$= \frac{2}{N(N-1)} \left[ \sum_{j=1}^{N-1} \alpha_j^{(1)} + \sum_{j=2}^{N-1} \alpha_j^{(1)} + \cdots + \sum_{j=N-1}^{N-1} \alpha_j^{(1)} \right]$

$$\leq \frac{2}{N(N-1)} \left[ \sum_{j=1}^{N-1} \alpha_j^{(2)} + \sum_{j=2}^{N-1} \alpha_j^{(2)} + \cdots + \sum_{j=N-1}^{N-1} \alpha_j^{(2)} \right]$$

$$= \frac{2}{N(N-1)} \sum_{j=1}^{N-1} j\alpha_j^{(2)} = \overline{d_2}$$

Corollary. It follows from the previous theorem that the average distance between nodes in an N-node network is at most equal to the average distance in an N-node chain, namely $\frac{(N+1)}{3}$ (see the appendix for the simple calculation of this value).

Remark. If G(A) denotes the Gini index of the array A of distance frequencies, we have

$$G(A) = \frac{1}{N-1} \left( N - 2\,\bar{d} \right) \qquad (13)$$

Hence, the Gini coefficient respects the extended majorization order. From (13) one can express $\bar{d}$ as a function of G(A):

$$\bar{d} = \frac{N}{2} - G(A) \frac{N-1}{2} \qquad (14)$$

The previous theorem shows that the operation of taking the average distance in an N-node network respects the opposite of the Lorenz majorization order, while the Gini coefficient respects this order.

**The median distance and its relation with the average distance in a chain**

Assume that we have an N-node chain, hence containing N-1 links. Then its set of distances contains $\frac{(N-1)N}{2}$ numbers and the median, Md, is either a natural number m or m-0.5. Then we have

$$1 + 2 + \cdots + (N - m) \geq \frac{N(N-1)}{4} > 1 + 2 + \ldots + (N - m - 1).$$

As for each natural number j, we have $\sum_{k=1}^{N-j} k = \frac{(N-j)(N-j+1)}{2}$, we can prove that m = [x] with

$$\frac{(N - x)(N - x + 1)}{2} = \frac{N(N - 1)}{4}$$

from which it follows that $x = \frac{(2N+1) - \sqrt{2N^2 - 2N + 1}}{2}$ and hence Md is either [x]-05 or [x]. For N large this leads to

$$Md \approx N\left(1 - \frac{\sqrt{2}}{2}\right) \approx 0.293\, N$$

Consequently, $\lim_{N \to \infty} \frac{Md}{\bar{d}} = 3\left(1 - \frac{\sqrt{2}}{2}\right) \approx 0.879 < 1$.

Moreover, we see that $Md < \bar{d} \Leftrightarrow N + \frac{1}{2} - \frac{\sqrt{2}}{2}\sqrt{N^2 - N + \frac{1}{2}} < \frac{N+1}{3} \Leftrightarrow N^2 - 13N + 4 > 0 \Leftrightarrow N > 12.7$. Hence, in practice: N ≥ 13. Checking this manually for N = 2, …,14 we find that also then $Md < \bar{d}$ except for N = 2, 5, 8, and 11 in which cases $Md = \bar{d}$.

**Returning to the example $G_0$ shown in Fig.2**

The $\alpha$ −array of $G_0$ is (6,7,6,2,0,0). Fig. 3 shows its extended Lorenz curve, situated between the extended Lorenz curve of $K_7$ (the complete network on 7 nodes) and the extended Lorenz curve of the chain of length 6. The average distances are respectively equal to 1, 2.19, and 2.67; the medians are 1, 2, and 2; while the corresponding Gini coefficients are: 0.833, 0.437, and 0.278.

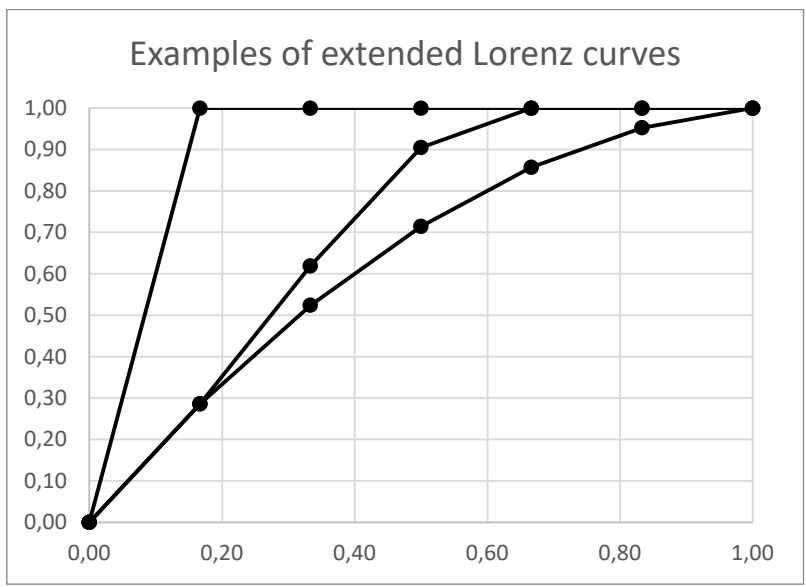

Fig. 3. Extended Lorenz curves of $K_7$, $G_0$ and $C_6$ (the chain of length 6)

**Conclusion**

In this article, we introduced the study of the distance distribution of a network. We showed that the distance distribution in an undirected network majorizes the one of a chain and is always smaller (in the sense of majorization) than the distribution of the corresponding complete N-network. The Gini coefficient respects the majorization order for such distributions, while the average distance behaves oppositely. As a consequence, the average and median distances in any such network are smaller than those of a chain.

We intend to use these results in the study of small worlds and the so-called six degrees of separation property (work in preparation).

Acknowledgement. The author thanks his colleague Ronald Rousseau for useful discussions.

**References.**

Appendix

Proof that the average distance in an N-node chain is $\frac{(N+1)}{3}$

The average distance in an N-node chain is equal to

$$\frac{2}{N(N-1)}\sum_{i=1}^{N-1} i(N-i) = \frac{2}{N(N-1)}\left(\frac{N^2(N-1)}{2} - \frac{(N-1)N(2N-1)}{6}\right)$$

$$= \frac{2N(N-1)}{2N(N-1)}\left(N - \frac{2N-1}{3}\right) = \frac{N+1}{3}.$$